\documentclass[conference]{IEEEtran}
\usepackage{amsfonts,latexsym,amsmath,amstext,amssymb,verbatim,epsfig,psfrag}
\usepackage{colordvi,pstcol}
\usepackage{graphicx}
\usepackage{psfrag,graphicx}

\newcommand\II{\mathcal{I}}

\def\R{\mbox{I\hspace{-.15em}R}}

\def\R{\mbox{I\hspace{-.15em}R}}

\def\Pb{\bold{P}}

\def\Qb{\bold{Q}}
\def\Ib{\bold{I}}
\def\Jb{\bold{J}}

\newtheorem{remark}{Remark}

\begin{document}
\title{Note: interpreting iterative methods convergence with diffusion point of view}
%\numberofauthors{1}
%\author{
%   \alignauthor Dohy Hong\\
%   \affaddr{Alcatel-Lucent Bell Labs}\\
%   \affaddr{Route de Villejust}\\
%   \affaddr{91620 Nozay, France}\\
%   \email{dohy.hong@alcatel-lucent.com}
%}
\author{\IEEEauthorblockN{Dohy Hong}
\IEEEauthorblockA{Alcatel-Lucent Bell Labs\\
Route de Villejust\\
91620 Nozay, France\\
dohy.hong@alcatel-lucent.com}
}

\date{\today}
\maketitle

\begin{abstract}
In this paper, we explain the convergence speed of different iteration schemes
with the fluid diffusion view when solving a linear fixed point problem.
This interpretation allows one to better understand why power iteration or Jacobi
iteration may converge faster or slower than Gauss-Seidel iteration.
\end{abstract}
%\category{G.1.0}{Mathematics of Computing}{Numerical Analysis}[Parallel algorithms]
%\category{G.1.3}{Mathematics of Computing}{Numerical Analysis}[Numerical Linear Algebra]
%\category{C.2.4}{Computer Systems Organization}{Computer-Communication Networks}[Distributed Systems]
%\terms{Algorithms, Performance}

%\keywords{Distributed computation, Iteration, Fixed point, Eigenvector.}

%plan

%%%%\section{Introduction}\label{sec:intro}
%%%%\section{Algorithm description}\label{sec:def}
%%%%\section{Experimental evaluation}\label{sec:eval}
%%%%\section{Conclusion}\label{sec:conclusion}

\begin{psfrags}
%%%%%%%%%%%%%%%%%%%%%%%%%%%%%%%%%%%%%%%%%%%%%%%%%%%%%%%%%%%%
\section{Introduction}\label{sec:intro}
%%%%%%%%%%%%%%%%%%%%%%%%%%%%%%%%%%%%%%%%%%%%%%%%%%%%%%%%%%%%
Based on the previous research results on the diffusion approach \cite{d-algo}
to solve fixed point problem in linear algebra, we propose here a new analysis of the
convergence speed of different iteration methods.

In Section \ref{sec:def}, we define the iteration methods that are considered. 
Section \ref{sec:anal} shows how to define the associated equivalent diffusion
iteration. Section \ref{sec:exm} shows few examples to illustrate the application.
%%%%%%%%%%%%%%%%%%%%%%%%%%%%%%%%%%%%%%%%%%%%%%%%%%%%%%%%%%%%
\section{Algorithms description}\label{sec:def}
%%%%%%%%%%%%%%%%%%%%%%%%%%%%%%%%%%%%%%%%%%%%%%%%%%%%%%%%%%%%

\subsection{Notations}

We will use the following notations:
\begin{itemize}
\item $\Pb \in \R^{N\times N}$ a real matrix;
\item $\Ib\in \R^{N\times N}$ the identity matrix;
\item $\Jb_i$ the matrix with all entries equal to zero except for
  the $i$-th diagonal term: $(\Jb_i)_{ii} = 1$;
\item $\Omega = \{1,..,N\}$;
\item $\sigma : \R^N \to \R$ defined by $\sigma(X) = \sum_{i=1}^N x_i$;
\item $\II = \{i_1,i_2,..,i_n,...\}$ a sequence of coordinate: $i_k \in \Omega$;
\item a fair sequence is a sequence where all elements of $\Omega$ appears infinitely often;
\item $e = (1/N,..,1/N)^T$.
\end{itemize}

\subsection{Problem to solve}
We will consider two types of linear fixed point problems:
$$
X = \Pb X
$$
$X\in \R^N$ and:
$$
X = \Pb X + B 
$$
$B\in \R^N$.

\subsection{Linear equation: $X = \Pb X$}
\subsubsection{Power iteration (PI)}
The power iteration $PI(\Pb, X_0)$ is defined by:
\begin{align}
X_{n} &= \Pb X_{n-1}\label{eq:pi}
\end{align}
starting from $X_0$.

\subsubsection{Gauss-Seidel iteration (GSl)}
Given a sequence of nodes for the update $\II = \{i_1,i_2,..,i_n,...\}$,
the Gauss-Seidel iteration $GSl(\Pb, X_0, \II)$ is defined by:
\begin{align}
(X)_{i_n} &= (\Pb X)_{i_n}\label{eq:gs-l}
\end{align}
starting from $X_0$, which simply means that the $n$-th update on $X$ is
on coordinate $i_n$ based on the last vector $X$ (each update modifying only
one coordinate of $X$). We could equivalently write it as:

\begin{align*}
X_n &= X_{n-1} + \Jb_{i_n} (\Pb - \Ib)  X_{n-1}.
\end{align*}

%%%
\subsection{Affine equation: $X = \Pb X + B$}
\subsubsection{Jacobi iteration (Jac)}
The Jacobi iteration $J(\Pb, B)$ is defined by:
\begin{align}
X_{n} &= \Pb X_{n-1} + B\label{eq:jacobi}
\end{align}
starting from $X_0 = (0,..,0)^T$.

\subsubsection{Gauss-Seidel iteration (GSa)}
Given a sequence of nodes for the update $\II = \{i_1,i_2,..,i_n,...\}$,
the Gauss-Seidel iteration $GSa(\Pb, B, \II)$ is defined by:
\begin{align}
(X)_{i_n} &= (\Pb X_{n-1} + B)_{i_n}\label{eq:gs-a}
\end{align}
starting from $X_0 = (0,..,0)^T$.
We could equivalently write it as:

\begin{align*}
X_n &= (\Ib - \Jb_{i_n}) X_{n-1} + \Jb_{i_n} (\Pb X_{n-1} + B).
\end{align*}

Note that in the above notations, the $n$ of $X_n$ may mean the
$n$-th application of an operator on a vector (Jacobi style, vector level update or VLU) or the $n$-th update
of a coordinate of $X$ (Gauss-Seidel style, coordinate level update or CLU).
In terms of diffusion view, the VLU can be always interpreted as a partial diffusion of CLU approach
cf. \cite{note}.

\begin{remark}
One could also consider a more general iteration scheme such as BiCGSTAB or GMRES,
but they require a fluid injection method, which may be more complex.
This is let for a future research.
\end{remark}
%%%%%%%%%%%%%%%%%%%%%%%%%%%%%%%%%%%%%%%%%%%%%%%%%%%%%%%%%%%%
\section{Diffusion equations}\label{sec:anal}
%%%%%%%%%%%%%%%%%%%%%%%%%%%%%%%%%%%%%%%%%%%%%%%%%%%%%%%%%%%%
The diffusion equations are defined by two state vectors $F_n$ and $H_n$
associated to the affine equation $X = \Pb X + B$ ($DI(\Pb, B, \II)$):

\begin{eqnarray}
F_n &=& F_{n-1} + (\Pb - \Ib) \Jb_{i_n} F_{n-1},\label{eq:diff}\\
H_n &=& H_{n-1} + \Jb_{i_n} F_{n-1}.\nonumber
\end{eqnarray}

The VLU adaptation of the above equation would be:
\begin{eqnarray}
F_n &=& \Pb F_{n-1},\label{eq:diff-vlu}\\
H_n &=& H_{n-1} + F_{n-1}.\nonumber
\end{eqnarray}
which is equivalent to
\begin{eqnarray*}
H_n &=& \Pb H_{n-1} + B,
\end{eqnarray*}
which is a Jacobi iteration.

\subsection{Linear equation: $X = \Pb X$}
\subsubsection{Power iteration (PI)}
We define: $F_0 = \Pb X_0 - X_0$, $H_0 = 0$ and the iterative diffusion equation
with VLU:
\begin{eqnarray*}
F_n &=& \Pb F_{n-1},\\
H_n &=& H_{n-1} + F_{n-1}.
\end{eqnarray*}
Then, we have the equalities:

\begin{eqnarray*}
F_n &=& X_{n+1} - X_n,\\
H_n &=& X_n - X_0.
\end{eqnarray*}

\subsubsection{Gauss-Seidel iteration (GSl)}
We set $F_0 = \Pb X_0 - X_0$, $H_0 = 0$ and $H_n, F_n$ defined by Equations \eqref{eq:diff},
then we have $X_n = H_n + X_0$.

\subsection{Affine equation: $X = \Pb X + B$}
\subsubsection{Jacobi iteration (Jac)}
The Jacobi iteration is equivalent to the VLU of diffusion equations \ref{eq:diff-vlu}.

\subsubsection{Gauss-Seidel iteration (GSa)}
The Gauss-Seidel iteration is equivalent to the (CLU of) diffusion equations \ref{eq:diff}.

\subsection{Summary of results on the equivalent iterations}
\begin{table}
\begin{center}
\begin{tabular}{|l||c|c|c|c|}
\hline
Type      &  VLU & CLU & VLU & CLU\\
\hline
Scheme    & PI                 & GSl         & Jac    & GSa\\
\hline
\hline
Defined by          & $\Pb, X_0$           & $\Pb, X_0,\II$ & $\Pb, B$ & $\Pb, B, \II$\\
\hline
Diffusion & $\Pb, \Pb X_0 - X_0$ & $\Pb, \Pb X_0 - X_0, \II$ & $\Pb, B$  & $\Pb, B, \II$.\\
\hline
\end{tabular}\caption{Iteration scheme equivalence.}\label{tab:equi}
\end{center}
\end{table}

%%% Exm
\section{Example of convergence speed comparison}\label{sec:exm}

For comparison analysis, we consider a PageRank equation \cite{page}:
$$
X = \Pb X = (d\Qb + (1-d)/N\Jb) X = d \Qb X + (1-d) e
$$
assuming $\sigma(X) = 1$.
Then we have:
\begin{itemize}
\item $PI(\Pb, X_0) = DI(d\Qb, d\Qb X_0 +(1-d) e-X_0)+X_0$;
\item $Jac(d \Qb, (1-d)e) = DI(d \Qb, (1-d)e)$;
\item $GSa(d \Qb, (1-d)e, \II) = DI(d \Qb, (1-d)e, \II)$;
\item $GSl(\Pb, X_0=e, \II) = DI(d \Qb, d(\Qb e-e), \II)$.
\end{itemize}

Note that $GSl(\Pb, X_0, \II)$ has no reason to converge in general. Il will converge
if $\II$ is a negative or positive fair sequence cf. \cite{cv}. Here, the decomposition
of $\Pb$ guarantees that it always converges for any fair sequence $\II$.
Note also that $Jac(d \Qb, (1-d)e)$ and $GSa(d \Qb, (1-d)e, \II)$ define
a non-decreasing vector $X_n$ (positive fluid diffusion).

Below, we compare PI, Jac and GS in simple scenarios.
\subsection{Case 1}
Note that the convergence is below measured in residual fluid using the $L_1$ norm
of $F_n$ of the equivalent diffusion equation. For an easy comparison, the x-axis shows
the number of iterations: for VLU approaches, it is exactly the index $n$ of $F_n$, for
CLU approaches, the $L_1$ norm of $F_{5\times n}$ is shown. 
We take $d=0.85$ and
$$
Q = 
\begin{pmatrix}
0&	0&	0&	0&	0.5\\
1&	0&	0&	0&	0.5\\
0&	1&	0&	0&	0\\
0&	0&	1&	0&	0\\
0&	0&	0&	1&	0
\end{pmatrix}
$$

\begin{figure}[htbp]
\centering
\includegraphics[angle=-90, width=8cm]{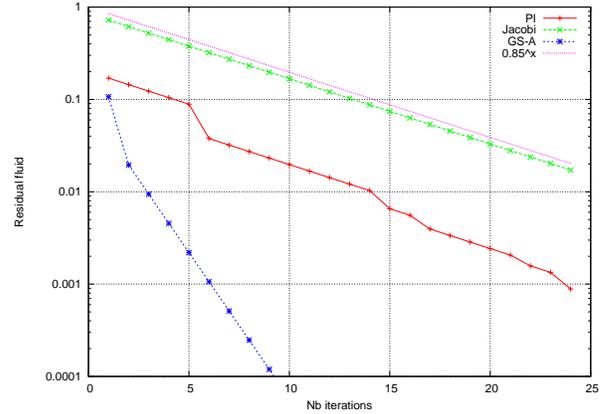}
\caption{Convergence.}
\label{fig:case1}
\end{figure}

Figure \ref{fig:case1}: Jacobi iteration converges exactly as $d^n$.
PI does roughly the same. GSa is much faster because we follow the graph
path and send cumulated fluid: it is roughly 4-5 times faster as expected.
GSl has the convergence slope of GSa but with a larger jump at the first iteration.

\subsection{Case 2}
We take $d=0.85$ and
$$
Q = 
\begin{pmatrix}
0&	0&	0.5&	0&	0.5\\
1&	0&	0&	0&	0.5\\
0&	1&	0&	0&	0\\
0&	0&	0.5&	0&	0\\
0&	0&	0&	1&	0
\end{pmatrix}
$$

\begin{figure}[htbp]
\centering
\includegraphics[angle=-90, width=8cm]{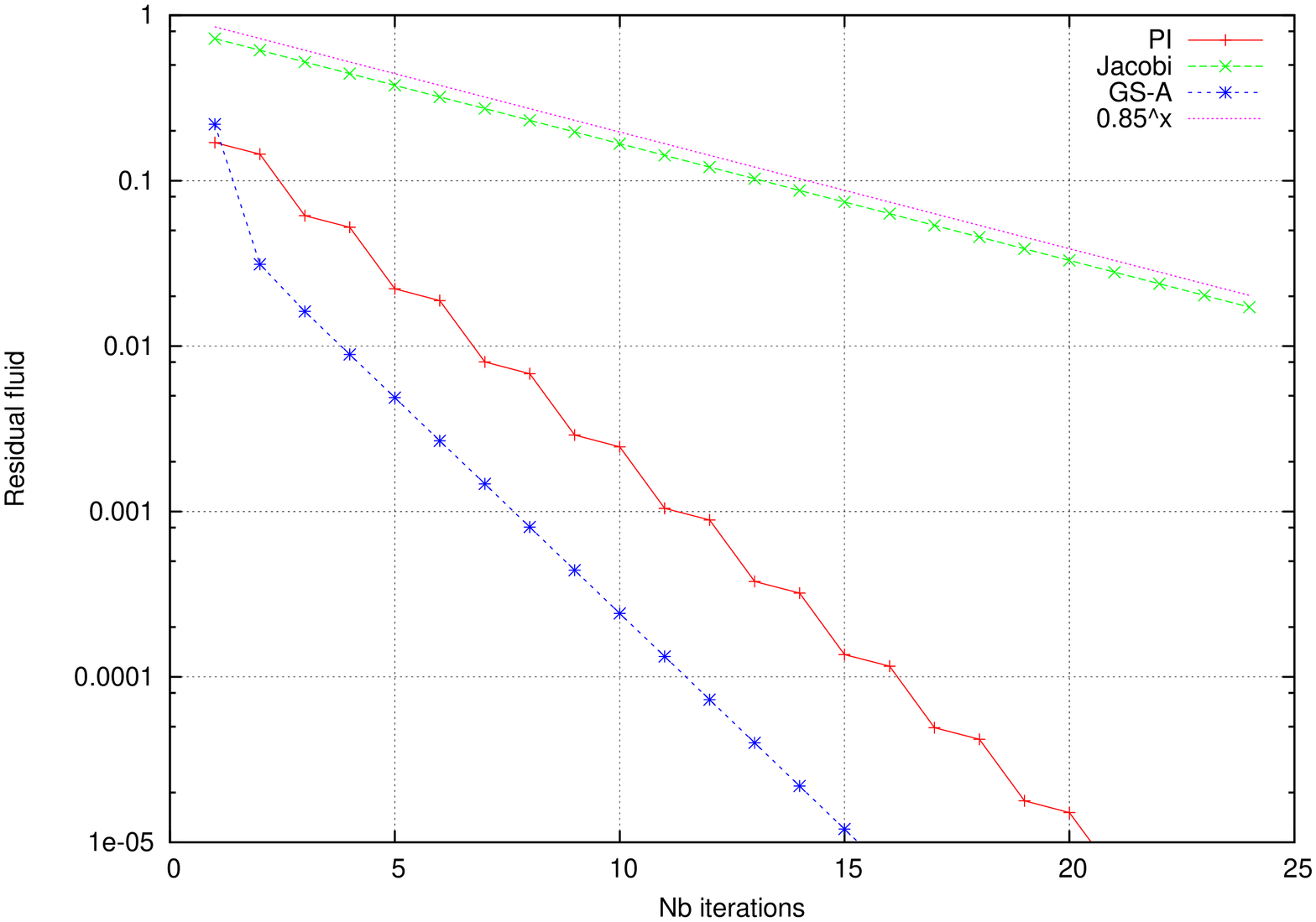}
\caption{Convergence.}
\label{fig:case2}
\end{figure}

Figure \ref{fig:case2}: Jacobi iteration converges still as $d^n$.
PI and GSa have similar convergence slope. GSl still shows the first jump.

\subsection{Case 3}
We take $d=0.85$ and
$$
Q = 
\begin{pmatrix}
0&	0.5&	0.5&	0&	0.5\\
1&	0&	0&	0&	0.5\\
0&	0.5&	0&	0&	0\\
0&	0&	0.5&	0&	0\\
0&	0&	0&	1&	0
\end{pmatrix}
$$

\begin{figure}[htbp]
\centering
\includegraphics[angle=-90, width=8cm]{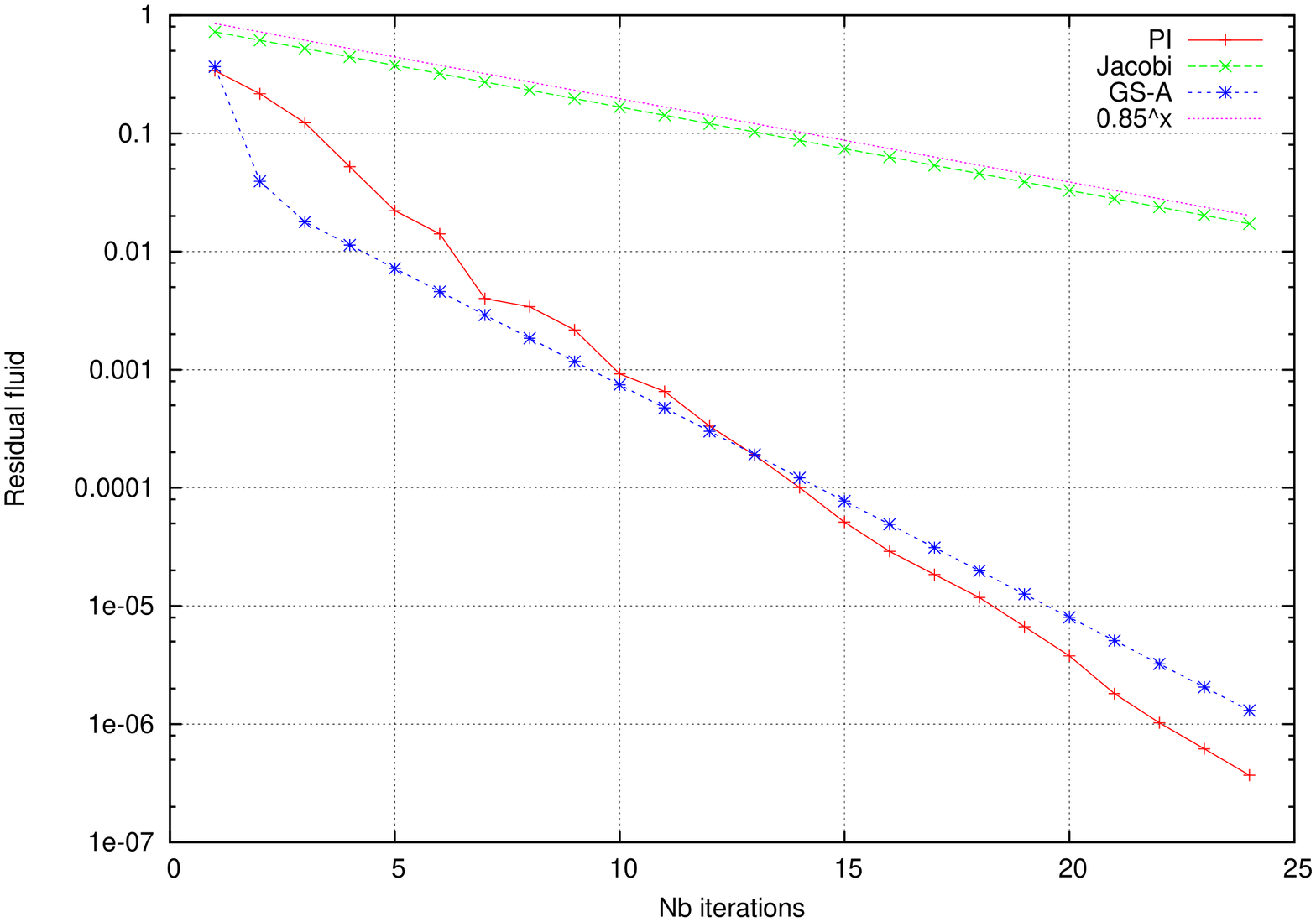}
\caption{Convergence.}
\label{fig:case3}
\end{figure}

Figure \ref{fig:case3}: Jacobi iteration converges always as $d^n$.
PI starts to converge faster than GSa when we add more links.
GSl has the first jump.

To explain, this convergence speed difference between the four methods, we
consider the case 4 below.

\subsection{Case 4}
We take $d=0.85$ and
$$
Q = 
\begin{pmatrix}
0&	0&	0&	0&	0.01\\
1&	0&	0&	0&	0\\
0&	1&	0&	0&	0\\
0&	0&	1&	0&	0\\
0&	0&	0&	1&	0.99
\end{pmatrix}
$$

\begin{figure}[htbp]
\centering
\includegraphics[angle=-90, width=8cm]{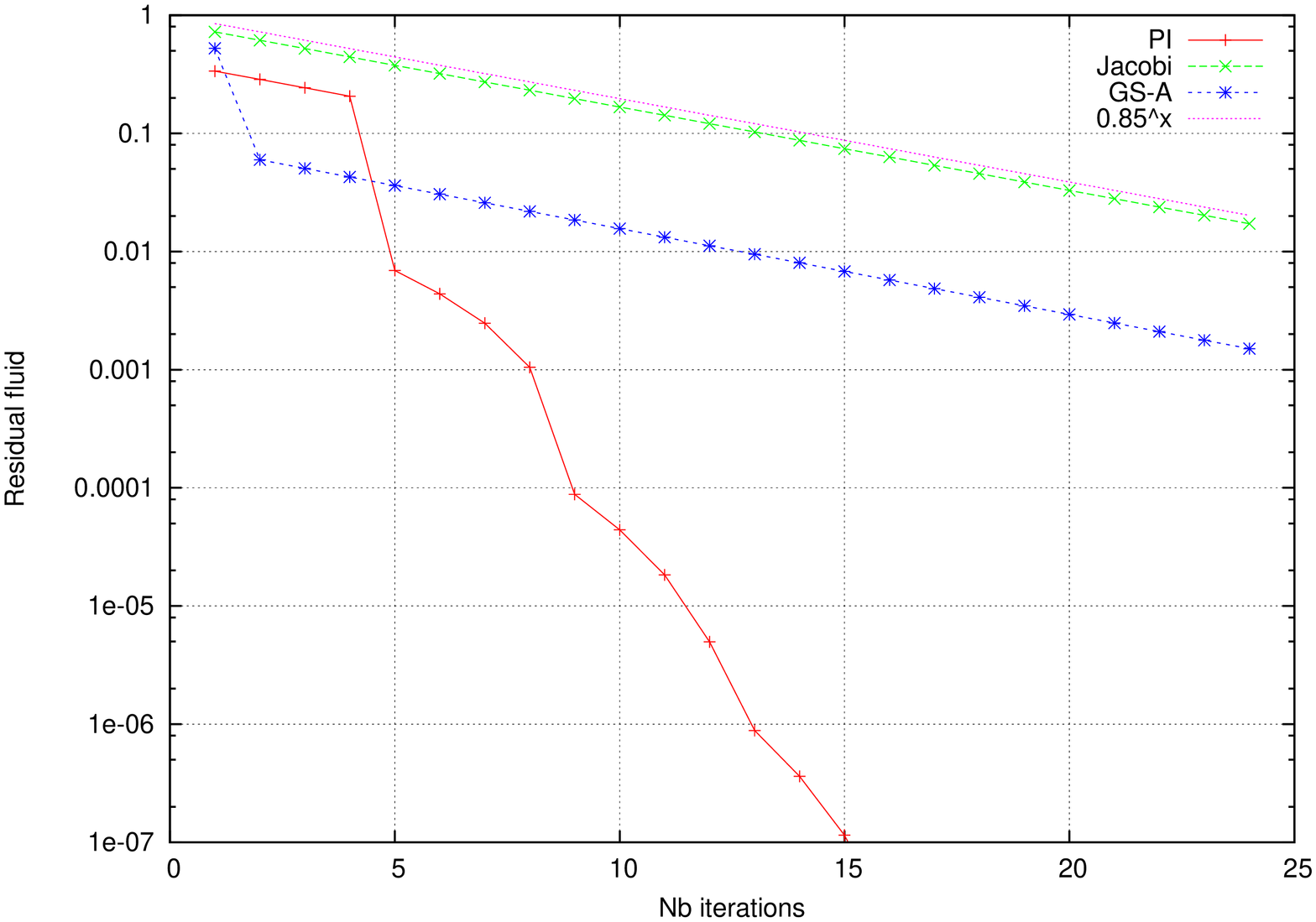}
\caption{Convergence.}
\label{fig:case4}
\end{figure}

Figure \ref{fig:case4}: with GSl, after one iteration (5 updates), we have a 
big jump due to the cumulated negative fluid that meets the positive fluid (node 5)
at 5-th update. Then, from 6-th update, we only move positive fluid. And only a fraction
of 1\% at node 5 are moved: this explains the convergence at $d^n$ after iteration 1 for
GSa and GSl. For cases 1-3, the whole fluid (no self loop) from each node is moved
to children nodes and this explains the gain factor (merging fluid before moving).

Now the good performance of PI can be explained by the fact that doing partial diffusion
we create more fluid cancellation and make the convergence faster than $d^n$.
To confirm this explanation, we plot in Figure \ref{fig:case4b} the amount of fluid that has been
canceled at each iteration. We see that this phenomenon is driving the convergence speed
(the difference to the residual fluid is due to the contracting factor that eliminates at each
diffusion a fraction $1-d$ (15\% here). In this case, the fluid disappeared after 5-th update
84\% due to fluid cancellation and 16\% due to the contracting factor.

\begin{figure}[htbp]
\centering
\includegraphics[angle=-90, width=8cm]{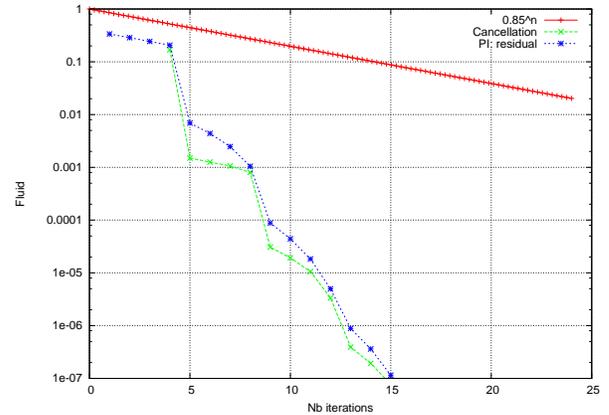}
\caption{Convergence: cancelled fluid.}
\label{fig:case4b}
\end{figure}
%%%%%%%%%%%%%%%%%%%%%%%%%%%%%%%%%%%%%%%%%%%%%%%%%%%%%%%%%%%%
\section{Conclusion}\label{sec:conclusion}
%%%%%%%%%%%%%%%%%%%%%%%%%%%%%%%%%%%%%%%%%%%%%%%%%%%%%%%%%%%%
In this paper, we described the equivalent equations of the diffusion iteration
associated to power iteration, Jacobi and Gauss-Seidel iteration and showed how
they can explain the convergence speed of each method.
\end{psfrags}
%======================================================================
\bibliographystyle{abbrv}
\bibliography{sigproc}

\end{document}